\documentclass[12pt]{article}
\usepackage{amsfonts}
\usepackage{amsmath}
\usepackage{amsthm}
\usepackage{amssymb}
\newcommand{\Complex}{\mathbb{C}}
\newcommand{\ComplexE}{\overline{\mathbb{C}}}
\newcommand{\Real}{\mathbb{R}}

\newcommand{\Int}{\mathbb{Z}}

\newcommand{\UD}{\mathbb{D}}
\newcommand{\UDe}{{\UD^*}}
\newcommand{\UTs}{{S^1}}
\newcommand{\UT}{\UTs}

\newcommand{\clS}{\mathcal{S}}
\newcommand{\clSt}{\clS^\infty}
\newcommand{\clSH}{\clS^{1,\alpha}}
\newcommand{\clSHo}[1]{\clS^{#1,\alpha}}
\newcommand{\clSqs}{\clS^{\mathrm{qc}}}

\newcommand{\mIm}{\mathop{\mathrm{Im}}}

\newcommand{\rotS}{\mathrm{Rot}(\UTs)}
\newcommand{\Kir}{\mathcal{M}}
\newcommand{\Mob}{\mathrm{M\ddot{o}b}(\UTs)}
\newcommand{\KirM}{K}
\newcommand{\DiffS}{\mathrm{Diff^+}(\UTs)}

\newcommand{\QS}{\mathrm{Homeo^+_{qs}}(\UTs)}
\newcommand{\id}{\mathrm{id}}

\newcommand{\T}[2]{\mathrm{T}_{#1}{#2}}

\newcommand{\der}{^\#}

\newcommand{\Lip}{\mathrm{Lip}_\alpha}
\newcommand{\Cc}{C^{1,\alpha}}
\newcommand{\Cco}[1]{C^{#1,\alpha}}
\newcommand{\Hol}{\mathrm{Hol}}
\newcommand{\HolC}{\mathrm{Hol}_C}

\newcommand{\Fr}{\mathcal F}
\newcommand{\KirJ}{J_\gamma}

\newcommand{\hder}{\partial}
\newcommand{\ader}{\bar\partial}

\newcommand{\sign}{\mathop{\mathrm{sgn}}}
\newtheorem{theorem}{Theorem}
\newtheorem{proposition}{Proposition}
\theoremstyle{remark}
\newtheorem{remark}{Remark}

\newcommand{\conj}[1]{\overline{\vphantom{\vbox to1.2ex{~}}#1}}

\newcommand{\res}{\mathop{\mathrm{Res\,}}}

\title{Matching univalent functions and conformal welding}
\author{Erlend\,Grong, Pavel\,Gumenyuk, and Alexander\,Vasil'ev}

\begin{document}
\date{}

\maketitle
{\let\thefootnote=\relax\footnotetext{2000 {\it Mathematics Suject Classification}.
Primary 30C35. Secondary 17B68.}\footnotetext{This research is supported by  {\it ESF Networking Programme} ''Harmonic and Complex Analysis and its Applications''; and by the the {\it Research Council of Norway}, grant \#177355/V30; second author is also supported by   the {\it Russian Foundation for Basic Research} (grant \#07-01-00120).}}

\begin{abstract}
Given a conformal mapping $f$ of the unit disk $\mathbb D$ onto a simply connected domain $D$  in the complex plane bounded by a closed Jordan curve, we consider the problem of constructing a matching conformal mapping, i.e., the mapping of the exterior of the unit disk  $\mathbb D^*$ onto the exterior domain $D^*$ regarding to $D$. The answer is expressed in terms of a linear differential equation with
a driving term given as the kernel of an operator dependent on the original mapping $f$. Examples are provided.
This study is related to the problem of conformal welding and to representation of the Virasoro algebra
in the space of univalent functions.
\end{abstract}

\section*{Introduction}

One of the classical problems of complex analysis  resides in finding the conformal mapping between a
given simply connected hyperbolic domain $D$ on the Riemann sphere $\ComplexE$ and some canonical domain,
e.\,g., the unit disk $\UD:=\{z:|z|<1\}$ or its exterior~$\UDe:=\ComplexE\setminus\overline\UD$, where 
$\overline\UD$ means the closure of $\UD$. Despite 
the fact that the existence and essential uniqueness of the mapping is guaranteed by the Riemann mapping
theorem, only in some particular cases it can be found analytically in a more or less explicit form. In
the present paper we consider a special formulation of this problem, when the domain~$D$ is bounded by a closed Jordan curve and represented by
means of the conformal mapping of $\UDe$ onto the exterior $D^*$ of the domain $D$, $\infty \in D^*$.

If the boundary $\partial D$ is $C^{\infty}$ smooth, then this formulation is
closely connected to Kirillov's representation of the Lie-Fr\'echet group $\DiffS$ of all
orientation preserving $C^\infty$-diffeomorphisms of the unit circle~$\UTs$, and to representation of the
Virasoro algebra, which is a central extension by $\Complex$ of the complexified Lie algebra of vector fields on $S^1$.
Virasoro algebra is known to play an important role in non-linear equations, where the Virasoro algebra is intrinsically related to the KdV canonical structure (see, e.g., \cite{Faddeev, Gervais}),  and in Conformal Field Theory, where the Virasoro-Bott group appears  as the space of reparametrization of a closed string  (see, e.g., \cite{Polchinski}) .

Let $f$  be a conformal mapping of $\UD$ onto a Jordan domain $D$ and $\varphi$ a conformal mapping of
$\UDe$ onto a Jordan domain $D^*$. The functions $f$~and~$\varphi$ are said to be {\it matching}\, if $D$
and $D^*$ are complementary domains, i.\,e., ${D\cap D^*=}$\O \  and ${\partial D=\partial D^*}$.

A pair of matching functions $(f,\varphi)$, being continuously extended to $\UTs$, defines a homeomorphism
of~$\UTs$ given by the formula
\begin{equation}\label{welding}
 \gamma=f^{-1}\circ \varphi.
\end{equation}
Such a representation of  homeomorphisms of $\UTs$ is called the {\it conformal welding}.

Using M\"obius transformations we can always assume that
\begin{itemize}
\item[(i)] $0\in D$ and $\infty\in D^*$;
\item[(ii)] $f(0)=f'(1)-1=0$;
\item[(iii)] $\varphi(\infty)=\infty$.
\end{itemize}
Conformal weldings have close connection to theory of quasiconformal (q.\,c.) mappings. Denote by $\clS$ the class
of all univalent analytic functions $f$ in $\UD$ subject to condition~(ii), and let $\clSqs$ be the
subclass of~$\clS$ consisting of functions which can be extended to a quasiconformal homeomorphism
of~$\ComplexE$. If $f\in\clSqs$, then $\varphi$ also admits q.\,c.~extention to $\ComplexE$ and therefore
$\gamma\in\QS$, where $\QS$ stands for the group of all orientation preserving quasisymmetric (q.\,s.)
homeomorphisms of~$\UTs$, i.e., $\gamma$ satisfies
\begin{equation}\label{qs}
\sup\left\{\left|\frac{\gamma\big(e^{i(t+h)}\big)-\gamma\big(e^{it}\big)}{\gamma\big(e^{i(t-h)}\big)-\gamma\big(e^{it}\big)}\right|
:t,h\in\Real,\,0<|h|<\pi\right\}<+\infty.
\end{equation}
Moreover, it is known that  for any $\gamma\in\QS$ there exists a unique
conformal welding~\eqref{welding} under conditions~(i)--(iii). Given
$\gamma\in\QS$, the construction of the pair $(f,\varphi)$ of matching functions  involves solution of the Beltrami equation $$ \ader f=\mu\,\hder f, $$ where $\hder$ and
$\ader$ stand for $\big(\frac\partial{\partial x}\mp i\frac\partial{\partial y}\big)/2$ respectively,
with the coefficient ${\mu=\mu(z)}$ depending on $\gamma$. See Section~\ref{QSSect} for details.

Some further study of the existence and uniqueness of conformal welding can be found in~\cite{Jones}.

In this paper we establish a more explicit connection between~$f$, $\varphi$ and~$\gamma$. We will use the
notation   $\Lip$, $\alpha\in(0,1)$ for the class of  H\"older continuous functions  of exponent~$\alpha$,
and $\Cco{n}$ for the class of $n$-times differentiable functions with the $n$-th derivative from the class
$\Lip$.  In order to indicate the domain of definition and admissible values of  functions  we will add them
in the parenthesis, e.\,g.,  $\Lip(\UTs,\Real)$ will stand the set of all real-valued functions which are from the
class $\Lip$ on $\UTs$. By $\clSHo{n}$, $n\geqslant1$, we denote the class of all functions $f\in\clS$
that map $\UD$ onto domains bounded by $\Cco{n}$-smooth Jordan curves. According to the
Kellog\,--\,Warschawski theorem (see, e.g.~\cite[p.\,49]{Pommer}), $f\in\clSHo{n}$ if and only if it can be
continuously extended  to $\UTs$, with $f|_\UTs\in\Cco{n}$, and  $f'|_\UTs$ does not vanish. The class of all
$f\in\clS$ that map~$\UD$ onto domains bounded by $C^\infty$-smooth Jordan curves will be denoted by
$\clSt$.

Let $f\in\clSH$. Consider the linear operator $I_f$ from~$\Lip(\UTs,\Real)$
to the space $\Hol(\UD)$ of all holomorphic functions in~$\UD$, defined by the formula
\begin{equation}\label{operator}
I_f[v](z):=-\frac{1}{2\pi i}\int_\UT\left(\frac{ s f'( s)}{f( s)}\right)^2\frac{v( s)}{f( s)-f(z)}\frac{d s}{ s},\quad z\in\UD.
\end{equation}

The following statement is our main result.
\begin{theorem}\label{thrm_main}
Suppose  $f\in\clSH$ and $\varphi$, $\varphi(\infty)=\infty$, are matching univalent functions. Then the kernel of
the  operator ${I_f:\Lip(\UTs,\Real)\to\Hol(\UD)}$ is the one-dimensional manifold\, $\ker
I_f=\mathrm{span}\{v_0\}$, where
\begin{equation}\label{v0}
v_0(z):=\frac1z\,\frac{(\psi\circ f)(z)}{f'(z)(\psi'\circ f)(z)},~~~\psi:=\varphi^{-1},\quad z\in\UTs.
\end{equation}
Moreover, the function $v_0$ is positive on
$\UTs$ and satisfies the condition
\begin{equation}\label{v_norm}
\int_0^{2\pi}\frac{dt}{v_0(e^{it})}=2\pi.
\end{equation}
\end{theorem}
\begin{remark}
Let $f\in\clSH$ be given. Consider the problem of finding the conformal mapping $\psi$
 of $D^*:=\ComplexE\setminus f(\UD)$ onto $\UDe$, $\psi(\infty)=\infty$, (subject to an additional condition  ensuring  the
uniqueness). Theorem\,\ref{thrm_main} reduces this problem to solution of the equation $I_f[v]=0$. Indeed,
given $f$ and $v_0$, one can calculate $\psi$ on the boundary of~$D^*$ by solving the following
differential equation $$\psi'(u)=H(u)\psi(u),\quad u\in\partial D^*,$$ where $H:=\tilde H\circ f^{-1}$ and
$\tilde H(z):= 1/\left[zf'(z)v_0(z)\right]$, $z\in\UTs$.
\end{remark}

Theorem\,\ref{thrm_main} describes the real-valued  solutions to the equation $I_f[v]=0$. The set of
complex solutions to this equation is much more extensive. Denote by $\HolC(\UDe)$ the class of all
continuous functions $h:\UDe\cup\UTs\to\Complex$ which are analytic in $\UDe$.

\begin{theorem}\label{thrm_complex}
Suppose  $f\in\clSH$ and $\varphi$, $\varphi(\infty)=\infty$, are matching univalent functions, and
$\gamma:=f^{-1}\circ\varphi$ is  the induced homeomorphism of $\UTs$. Then the kernel of the operator
${I_f:\Lip(\UTs,\Complex)\to\Hol(\UD)}$ coincides with the set of all functions $v$ of the form
\begin{equation}\label{vComplex}
 v(z)=v_0(z)\cdot(h\circ \gamma^{-1})(z),\quad z\in\UTs,
\end{equation}
where $h$ is an arbitrary function belonging to $\HolC(\UDe)\cap\Lip(\UTs,\Complex)$ and $v_0$ is defined by~\eqref{v0}.
\end{theorem}

In Section\,\ref{KirSect} we show how the operator $I_f$ appears in a natural way within the
identification of the Kirillov's homogeneous manifold $\Kir:=\DiffS/\rotS$ with $\clSt$ and deduce an
analogue of Theorem\,\ref{thrm_main} for the $C^\infty$-smooth case.

Section\,\ref{proof} is devoted to the proof of Theorems\,\ref{thrm_main}\,and\,\ref{thrm_complex}.
Examples of univalent matching functions and conformal weldings are given in Sections~5 and~6.

\section{Conformal welding for quasisymmetric homeomorphisms of {\mathversion{bold}$\UTs$}}\label{QSSect}

It is known that {\it conformal welding establishes a bijective correspondence between $\clSqs$ and
$\QS/\rotS$},  where $\rotS$ stands for the group of rotations of~$\UTs$. For the history of the question,
see e.\,g.~\cite{GFT}. Here we briefly give a sketch of the proof, see also \cite{Takh}.

Let $u$, $u(\infty)=\infty$, be any q.\,c.~automorphism of $\UDe$. Let us construct the quasiconformal
homeomorphism  $\tilde f$ of the Riemann sphere~$\ComplexE$, such that the functions $f:=\tilde f|_\UD$ and
$\varphi:=\big(\tilde f|_\UDe\big)\circ u$ are analytic in $\UD$ and $\UDe$ respectively. It is easy to
see that $\tilde f$ should satisfy the Beltrami equation
\begin{equation}\label{Beltrami}
\ader \tilde f(z)=\mu(z)\,\hder \tilde f(z),\quad \mu(z):=\left\{\begin{array}{ll}\ader
\big(u^{-1}(z)\big)/\hder\big(u^{-1}(z)\big), & \text{if $z\in\UDe$,}\\0, &
\text{otherwise.}\end{array}\right.
\end{equation}

In order to have a unique solution we impose the following normalization
\begin{equation}\label{normirovka}
\tilde f(0)=\tilde f'(0)-1=0,~~~ \tilde f(\infty)=\infty.
\end{equation}
Then $f\in\clSqs$ and $\varphi$ are
matching functions and the homeomorphism of the unit circle $\gamma:=f^{-1}\circ\varphi$ coincides with
the continuous extension of $u$ to~$\UTs$.

It is known~\cite{BeurAhlf} that an orientation preserving homeomorphism ${\gamma:\UTs\to\UTs}$ can be
extended to a q.\,c.~automorphism $u$ of $\UDe$ if and only if it is quasisymmetric, i.\,e.,
satisfies~\eqref{qs}. Moreover, by  superposing  $u$ and a suitable q.\,c. automorphism  of~$\UDe$, 
identical on~$\UTs$,  one can always assume that $u(\infty)=\infty$. It follows that for any $\gamma\in\QS$
there exists a  conformal welding with $f\in\clSqs$.

Fix any q.\,c. extension $u:\UDe\to\UDe$; $\infty\mapsto\infty$, of $\gamma\in\QS$ and let $$\tilde
f(z):=\left\{\begin{array}{ll}f(z),&\text{if $z\in\UD$},\\(\varphi\circ
u^{-1})(z),&\text{otherwise,}\end{array}\right.$$ where $f\in\clS$ and $\varphi$ are matching univalent
functions such that $\gamma=f^{-1}\circ\gamma$. Then $\tilde f$
satisfies~\eqref{Beltrami}\,--\,\eqref{normirovka}. This defines $\tilde f$ uniquely~(see, e.\,g.,~\cite[p.\,194]{LehtoBook}). It follows that for any $\gamma\in\QS$ the conformal welding is unique.

On the hand, if $f\in\clSqs$, then $\varphi$ and consequently $\gamma=f^{-1}\circ\varphi$, can be extended
to a quasiconformal homeomorphism of~$\ComplexE$. It follows that ${\gamma\in\QS}$. Since the condition
$\phi(\infty)=\infty$ defines a conformal mapping onto $D^*:=\ComplexE\setminus\overline{f(\UD)}$ only up to
rotations, $f$ corresponds to the equivalence class $[\gamma]\in\QS/\rotS$, rather than an element of
$\QS$.
\begin{remark}
If $\gamma:\UTs\to\UTs$ is a diffeomorphism, then one of its q.\,c. extensions  to $\UDe$ is given by
the formula $u(r e^{it}):=r\gamma(e^{it})$, and the Beltrami coefficient $\mu$ in~\eqref{Beltrami} equals
$$\mu(re^{it})=e^{2it}\,\frac{1-\big(\gamma^{-1}\big)\der(e^{it})}{1+\big(\gamma^{-1}\big)\der(e^{it})},$$
where we introduce an operator `$\#$' by  $\beta^\#:=\big(\pi^{-1}\circ\beta\circ\pi\big)',$ and $\pi:\Real\to\UTs$ is the universal covering,
$\pi(x)=e^{i x}$.

\end{remark}

In Section~\ref{Poly} we consider a certain class of analytic diffeomorphisms~$\gamma$ for which
Theorem~\ref{thrm_main}  can be used to find the conformal welding without solving the Beltrami equation.

\section{Kirillov's representation of {\mathversion{bold}$\DiffS$} via univalent functions}\label{KirSect}

The group $\DiffS$ of all orientation preserving $C^\infty$-diffeomorphisms of the unit circle~$\UTs$ is
one of the simplest, and by this reason important, example of an infinite-dimensional Lie group.  Denote
by~$\Fr$ the Fr\'echet space  of all {$C^\infty$-smooth}  functions $h:\UTs\to\Real$ endowed with the
countable family of seminorms $\|h\|_n:=\max_{x\in\Real}\big|(d^n/dx^n)h(e^{ix})\big|$, $n\geqslant0$. It
is known (see,~e.\,g.,~\cite{LieFre}) that $\DiffS$ becomes a Lie-Fr\'echet group if we define  the structure of a $C^\infty$-smooth manifold on
$\DiffS$ by means of  the covering mapping $h\mapsto
\gamma[h]$, ${\gamma[h](\zeta):=\zeta e^{ih(\zeta)}}$, of the open set $\{h\in\Fr:dh(e^{ix})/dx>-1\}$ onto
$\DiffS$. All the tangent spaces $\T\gamma\DiffS$ are identified then in a natural way with~$\Fr$.

Kirillov~\cite{Kir} suggested to use the correspondence between $\QS$ and $\clSqs$ established by
means of conformal  welding, in order to represent the homogenous manifold $\Kir:=\DiffS/\rotS$,  usually
referred to as Kirillov's manifold, via univalent functions.

Consider the class~$\clSt$ of all functions $f\in\clS$  having $C^\infty$-smooth extension
to~$\partial\UD$ with non-vanishing derivative. By the Kellog\,--\,Warschawski theorem (see,
e.g.,~\cite[p.\,49]{Pommer}), $f\in\clSt$ if and only if $f$ has  a $C^\infty$-smooth extension to $S^1$ and the derivative $f'\big|_{S^1}$ does not vanish. It follows that $\clSt$ corresponds via conformal welding to a subset of $\DiffS/\rotS$. According
to the result of Kirillov~\cite{Kir}, it actually coincides with~$\DiffS/\rotS$, and consequently one can
identify $\Kir$ with $\clSt$.

Denote by~$\KirM:\clSt\to\Kir$ the mapping that takes  each $f\in\clSt$ to the corresponding equivalence
class of diffeomorphisms~$[\gamma]$. The infinitesimal version of the inverse mapping is as follows.

Fix any $v\in\Fr\cong\T\id\DiffS$ and consider the right-invariant vector field over $\DiffS$,
$V:\gamma\mapsto v\circ\gamma\in\Fr\cong\T\gamma\DiffS$ generated by $v$. This gives us the identification
$\T\gamma\DiffS\cong \T\id\DiffS\cong\Fr$, which we adhere further on,  and which is obviously
different from the identification of $\T\gamma\DiffS$ with $\Fr$ described above.

Thus, to each $v\in\Fr$ and each ${\gamma\in\DiffS}$ one associates the
variation~$\gamma_\varepsilon(\zeta):=\gamma(\zeta)\exp[i\varepsilon (v\circ\gamma)(\zeta)]$ of~$\gamma$.
According to~\cite{Kir1}, the corresponding variation of the function~$f$ equals to
$f_\varepsilon:=\KirM^{-1}([\gamma_\varepsilon])=f+\delta f+o(\varepsilon)$, where
\begin{equation}\label{deltaf}
\delta f(z)=\frac{\varepsilon}{2\pi }\int_\UT\left(\frac{ s f'( s)}{f( s)}\right)^2\frac{
f^2(z)\,v(s)}{f(z)-f(s)}\frac{d s}{ s}=i\varepsilon f^2(z) I_f[v](z),\quad z\in\UD.
\end{equation}

A natural consequence  is that $I_f[v](z)=0$ for all $z\in\UD$ if and only if the variation of
$\gamma$  produces no variation of $[\gamma]\in\Kir$ (up to higher order terms). It can be reformulated
as follows: {\it  the element of~$\T\gamma\DiffS$ represented by $v\circ \gamma$ is tangent to the
one-dimensional manifold $$\gamma\circ\rotS=[\gamma]\subset\DiffS.$$} The latter is equivalent to
$$v\in\mathrm{Ad}_{\gamma}\Big(\T\id\rotS\Big)=\mathrm{Ad}_{\gamma}\big\{\text{constant functions on
$\UTs$}\big\}.$$ Elementary calculations show that $$\mathrm{Ad}_\gamma
u=\frac{u\circ\gamma^{-1}}{\big(\gamma^{-1}\big)^\#}.
$$
As a conclusion we get
\begin{proposition}\label{prop}
 The kernel of $I_f:\Fr\to\Hol(\UD)$ is one-dimensional and coincides with $\mathrm{span}\{1/(\gamma^{-1})^\#\}$.
\end{proposition}

\begin{remark}
Proposition~\ref{prop} reveals a version of Theorem~\ref{thrm_main} for $C^\infty$-smooth case. It
reduces the problem of calculating $\KirM^{-1}(f)$ to solution of the equation $I_f[v]=0$. The nontrivial
solution $v_0$ subject to the normalization $$\int_0^{2\pi}\frac{dt}{v_0(e^{it})}=2\pi$$ allows us to
determine $[\gamma]$ by means of the equality
\begin{equation}\label{gamma} \gamma^{-1}(e^{ix})=\exp
\left(\int_0^x\frac{i\,dt}{v_0(e^{it})}+iC\right),
\end{equation}
with the arbitrary constant $C$ being responsible for the fact that~\eqref{gamma} defines $\gamma$  only
up to the right action of $\rotS$.
\end{remark}

\section{Virasoro algebra and complex structure on Kirillov's manifold}
The Lie algebra of $\DiffS$ is the Fr\'echet space $\Fr$ endowed with the Lie bracket
\begin{equation}\label{braket}
 \{v_1,v_2\}(e^{ix})=v_2(e^{ix})\,\frac{dv_1(e^{ix})}{dx}-v_1(e^{ix})\,\frac{dv_2(e^{ix})}{dx}.
\end{equation}
\begin{remark}
The expression~\eqref{braket} differs in  sign from the commutator~$[V_1,V_2]$ of the vector fields
$V_j:\gamma\to v_j\circ\gamma$ generated by $v_j$, because $V_j$ are right-invariant vector fields rather
than left-invariant, which are usually considered in this context.
\end{remark}

The simplest basis for the complexification $\Fr_\Complex:=\{v_1+iv_2:v_1,v_2\in\Fr\}$ of $\Fr$ is 
given by powers of $z$: $$L_k(z):=iz^k,\quad k\in\Int.$$ Continuation of the Lie bracket
$\{\cdot,\cdot\}:\Fr\times\Fr\to\Fr$ by complex bilinearity to $\Fr_\Complex$ gives the commutation
relations $\{L_k,L_j\}=(j-k)L_{k+j}$.

The (complex) Virasoro algebra can defined now as the central extension of $\Fr_\Complex$ by $\Complex$
which  is the Lie algebra over $\Fr_\Complex\oplus\Complex$ with the commutation relations
$$\big\{\big(L_k,a\big),\big(L_j,b\big)\big\}=\big(\{L_k,L_j\},\tfrac{c}{12} k(k^2-1)\delta_{k,-j}\big).$$
Here $c$ is a constant parameter referred to  as the {\it central charge} in Mathematical Physics.

Unfortunately, it is not known whether the Lie-Fr\'echet algebra $\Fr_\Complex$ is the Lie algebra of
any  Lie\,--\,Fr\'echet group, which, if exists, can serve as complexification for $\DiffS$. There are
strong reasons to believe that such a group does not exist~\cite{Neretin}. Nevertheless, the infinitesimal
action $\Fr\times\Kir\to\T{}\Kir$ induced by the left action of $\DiffS$, can be extended from $\Fr$ to
$\Fr_\Complex$, due to the fact that the linear space spanned by the variations~\eqref{deltaf} has a
natural complex structure, the operation of multiplication by $i$. This induces complex structure~$\KirJ$
on $\Fr/\ker I_f\cong\T{[\gamma]}\Kir$. We use Theorem\,\ref{thrm_complex} to obtain the explicit form of
it. Instead of looking for the operator on $\Fr/\ker I_f$ we define $\KirJ$ as an operator on $\Fr$ with
the property that $\KirJ[v_0]=0$. For $v\in\Fr$ we have $$iI_f[v]=I_f[iv]=I_f[\KirJ v].$$ It follows that
$\KirJ v=iv-\tilde v$, where $\tilde v\in\Fr_\Complex$ is a solution of $I_f[\tilde v]=0$ satisfying the
condition~$\mIm\tilde v=v$. Using the representation~\eqref{vComplex} for $\tilde v$ we obtain the formula
\begin{equation}
\KirJ[v]\circ\gamma=(v_0\circ\gamma)\cdot J_0\left[\frac{v\circ\gamma}{v_0\circ\gamma}\right],
\end{equation}
where $J_0:\Fr\to\Fr$ is the so-called conjugation, $$J_0\left[\sum_{k\in\Int}a_kz^k\right]=i\sum_{k\in\Int}\sign(k)a_kz^k.$$
Elementary calculations lead us to the following
\begin{proposition}
The complex structure on~$\T{}\Kir$ induced by the standard complex structure on $\Fr_\Complex$ via $I_f$
is given by $\KirJ=\mathrm{Ad}_\gamma J_0\left(\mathrm{Ad}_\gamma\right)^{-1}$, where $\mathrm{Ad}_\gamma$
stands for the differential of $A_\gamma\beta:=\gamma\circ\beta\circ\gamma^{-1}$ at the
origin~$\beta=\id$.
\end{proposition}
\begin{remark}
The complex structure~$\KirJ$ coincides with that introduced in~\cite{AirMall} only for the case
$\gamma=\id$ and thus it is not invariant under the right action of $\DiffS$ on $\Kir$. However, $\KirJ$
is left-invariant, which is proved by Kirillov~\cite{Kir} and easily follows from the fact that the
differential of the left action of $\DiffS$ is given by $v\mapsto \mathrm{Ad}_\gamma v$, where
$v\in\Fr\cong\T\gamma\DiffS$.
\end{remark}

\section{Proof of Theorems~\ref{thrm_main} and \ref{thrm_complex}}\label{proof}

Here we give a proof of Theorems~\ref{thrm_main} and \ref{thrm_complex} stated in the  Introduction, which
is based purely on complex analysis.

\begin{proof}[Proof of Theorem~\ref{thrm_main}]
Denote $D:=f(\UD)$, $\Gamma:=\partial D$, $$H(u):=\frac{g(u)v(g(u))}{u^2g'(u)},\quad F(w):=-\frac1{2\pi
i}\int_{\Gamma}\frac{H(u)}{u-w}du,~~  w\in \ComplexE\setminus\Gamma,$$ where $g$ stands for the inverse of
the function~$f$.

The equation $I_f[v](z)=0$, $z\in\UD$, is equivalent to
\begin{equation}\label{eq_eq}
 \quad F(w)=0,\quad w\in D.
\end{equation}
Using the Sokhotsky~--~Plemelj formulas we conclude that if $v$ is a solution to~\eqref{eq_eq}, then
$H(u)$ is the boundary values of an analytic function in $D^*:=\ComplexE\setminus\overline D$ vanishing at
$w=\infty$. The converse is also true due to the Cauchy integral formula for unbounded domains. It follows
that $v_0$ is a solution to~$\eqref{eq_eq}$. Indeed, for $v=v_0$ we have
$$H(u)=\frac{\psi(u)}{u^2\psi'(u)}.$$

The function~$v_0$ can be expressed as $v_0(z)=\zeta\varphi'(\zeta)/\big(zf'(z)\big)$, where
${\zeta:=\psi(f(z))}$. Both
vectors $\zeta\varphi'(\zeta)$ and $zf'(z)$ are the outer normal vectors of $\Gamma$ at the point
$w=f(z)=\varphi(\zeta)$. It follows that $v_0(z)>0$. The continuous function $\tau(t)$ defined by
${e^{i\tau(t)}=\psi\big(f(e^{it})\big)}$, $t\in\Real$,  satisfies the conditions $\tau'(t)=1/v_0(e^{it})$
and $\tau(t+2\pi)=\tau(t)+2\pi$. It follows that~\eqref{v_norm} holds.

It remains to prove that any real-valued solution~$v\in\Lip(S^1,\mathbb R)$ to equation~\eqref{eq_eq} is of  the form
$v=\lambda v_0$, $\lambda\in\Real$. Assume $v_1\in\Lip(S^1,\mathbb R)$ is a solution. And consider the one-parameter
family of solutions defined by $v:= v_0+\varepsilon v_1$, where $\varepsilon\in\Real$ is sufficiently
small for $v$ to be positive on $\UT$. By the above argument, the function $G(u):=u\psi'(u)H(u)$,
$u\in\Gamma$, has an analytic continuation to $D^*$, which will be denoted by $G(w)$.

The function $G$ does not vanish in $D^*\cup\Gamma$ provided $\varepsilon$ is small enough.  Indeed,
$G(w)\to\psi(w)/w$ as $\varepsilon\to 0$ uniformly in $D^*\cup\Gamma$, with the limit function $\psi(w)/w$
continuous and non-vanishing.  It follows that $\tilde G(w):=\log G(w)$ is analytic in $D^*$ and
continuous on $D^*\cup\Gamma$. The inequality $v>0$ implies that $\mIm \tilde G(u)=\mIm \log J(u)$,
${u\in\Gamma}$, where $J(u):=g(u)\psi'(u)/\big(u g'(u)\big)$. This equality  determines $\tilde G$ up to a
real constant term. Therefore, $v(z)$ is unique up to a positive constant coefficient. This completes the
proof.
\end{proof}

By the same techniques one can prove Theorem~\ref{thrm_complex}.
\begin{proof}[Proof of Theorem~\ref{thrm_complex}]
Let us look for solutions to $I_f[v]=0$ in the form~\eqref{vComplex}  {\it without any  a priori
assumptions} on $h$, except for that $h\in\Lip(\UTs,\Complex)$. Any solution can be represented in this
form because $v_0$ is positive. Now we use the change of variable $s=\gamma(t)$ in
integral~\eqref{operator}. Taking into account that
$$v_0(s)=1/(\gamma^{-1})\der(s)=(t/s)\cdot(ds/dt)\text{ and
}f'\big(\gamma(t)\big)\cdot(ds/dt)=\varphi'(t),$$
 we conclude that
$$ I_f[v](z)=-\frac{1}{2\pi i}\int_\UT\left(\frac{t
\varphi'(t)}{\varphi(t)}\right)^2\frac{h(t)}{\varphi(t)-w}\frac{d t}{t},\quad w:=f(z),~~ z\in\UD.$$

Applying another one change of variable $u=\varphi(t)$, we obtain the following expression for the above
quantity $$ -\frac{1}{2\pi i}\int_\Gamma\frac{\psi(u)}{u\psi'(u)}\,\frac{h(\psi(u))/u}{u-w}\,du , $$ Due
to the Sokhotsky-Plemelj formulas and the Cauchy integral formula for unbounded domains, the  above
quantity equals zero for all $w\in D:=f(\UD)$ if and only if $h$ represents the boundary values of an
analytic function in~$\UDe$. This fact proves the theorem.
\end{proof}

\section{Examples of matching univalent functions}

Here we consider a class of examples, for which both matching functions $f$ and $\varphi$ are expressed by
means of ordinary differential equations.

Given an integer $n>1$, let us consider the following quadratic differentials
\begin{align*}
&\Xi(\zeta) d\zeta^2:=-\frac{d\zeta^2}{\zeta^2};\\
&W(w)dw^2:=-\frac{w^{n-2}dw^2}{P(w)},\quad P(w):=\prod_{k=0}^{n-1}(w-w_k),\quad w_k:=e^{2\pi i
k/n};\\
&Z(z)dz^2:=-\frac{z^{n-2}dz^2}{Q(z)}; \quad Q(z):=\varkappa \prod_{k=0}^{n-1}\frac{|z_k|}{z_k}(z_k-z)(z-1/\conj{z_k}),\quad z_k:=r e^{2\pi i
k/n},\\
\end{align*}
where $r\in(0,1)$, and
 $\varkappa>0$ is such that $\int_\UTs\sqrt{Z(z)}dz=2\pi$ for the appropriately chosen branch of the
square root.

These quadratic differentials have the following structure of trajectories~(see e.\,g., \cite{Jenkins, AYu}). All the
trajectories of $\Xi(\zeta) d\zeta^2$ are circles centered on the origin, with $0$ and $\infty$ as
critical points. Critical trajectories of $W(w)dw^2$ are line intervals joining $w=0$ with $w_k$. Denote
the union of their closures by $E_w$. All the remaining trajectories are closed Jordan curves separating
$E_w$ and the critical point at infinity. The structure of trajectories of the quadratic differential
$Z(z)dz^2$ is symmetric with respect to the unit circle, which is also a trajectory. Similarly
to~$W(w)dw^2$, singular trajectories of~$Z(z)dz^2$ that lies in~$\UD$ are line intervals joining the
origin with~$z_k$. They form a continuum, which we denote by $E_z$. The singular trajectories lying
outside~$\UD$ form the symmetric continuum~$E_z^*$. All the remaining trajectories are Jordan curves
separating $E_z$ and $E_z^*$.

Let us choose any non-singular trajectory~$\Gamma$ of quadratic differential $W(w)dw^2$ and construct the
bijective conformal mappings $f:\UD\to D$, $f(0)=0$, $f'(0)>0$, and $\varphi:\UDe\to D^*$,
$\varphi(\infty)=\infty$, $\varphi'(\infty)>0$, where $D$ and $D^*$ are the interior and exterior of
$\Gamma$, respectively.

The mapping~$f$ can be constructed as follows. Let us define the parameter~$r$ in~$Z(z)dz^2$ by requiring
that the moduli of the annular domains $\UD\setminus E_z$ and $D\setminus E_w$ are equal. Consider the
conformal mapping $f$ of $\UD\setminus E_z$ onto $D\setminus E_w$ normalized by $f(z_0)=w_0$. This mapping
satisfies the following differential equation
\begin{equation}\label{diff_eq_f}
W(w)dw^2=Z(z)dz^2.
\end{equation}
Indeed, the conformal mapping~$\zeta=\varrho(z)$ of the ring domain $\ComplexE\setminus(E_z\cup E_z^*)$
onto the domain of the form $G:=\{\zeta:\rho<|\zeta|<1/\rho\}$ normalized by $\varrho(z_0)=\rho$ satisfies the
equation (see, e.\,g.~\cite[p.~43--46]{AYu})
\begin{equation}\label{diff_eq_rho}
Z(z)dz^2=\Xi(\zeta) d\zeta^2.
\end{equation}
Analogously, the conformal mapping $\zeta=\psi(w)$ of the circular domain ${\ComplexE\setminus E_w}$ of
the quadratic differential~$W(w)dw^2$ onto the domain $\{z:|z|>\rho\}$ normalized by $\psi(\infty)=\infty$
and $\psi(w_0)=\rho$ satisfies the equation
\begin{equation*}\label{diff_eq_psi}
W(w)dw^2=\Xi(\zeta) d\zeta^2.
\end{equation*}
Since the moduli of the annular domains $\UD\setminus E_z$ and $D\setminus E_w$ are equal,
$\psi(D\setminus E_w)=G'$, $G':=\{\zeta:\rho<|\zeta|<1\}$, and consequently $f=\psi^{-1}\circ\varrho$. It follows that~\eqref{diff_eq_f}
holds.

Now using the symmetry of $E_w$ and $E_z$ one can prove that $f$ extends analytically to $E_z$, i.e., $f$
is the desired conformal mapping of $\UD$ onto $D$.

It follows from the above consideration, that the exterior mapping is ${\varphi=\psi^{-1}|_\UDe.}$

By rescaling $w$-plane we can assure that $f\in \clS$. Now we can easily calculate the function $v_0$
spanning the kernel of the operator~$I_f[v_0]$, formula~\eqref{operator}. According to
Theorem~\ref{thrm_main} and equality~\eqref{diff_eq_rho},
$$v_0(z)=\left(-z^2Z(z)\right)^{-1/2}=\sqrt\frac{\varkappa}{r^n}\,\prod_{k=0}^{n-1}|z-re^{ikt/n}|,\quad z\in\UTs.$$
\begin{remark}
The choice of the coefficient~$\varkappa$ in the construction of quadratic differential~$Z(z)dz^2$ garantees that $v_0$ 
satisfies normalization~\eqref{v_norm}.
\end{remark}

\begin{remark}
The circle diffeomorphism~$\gamma$  coincides on~$\UTs$ with~$\varrho^{-1}$. Consequently, it can be
extended analytically from $\UTs$ to the ring~$G$.
\end{remark}

\begin{remark}
For the case $n=2$ the curve $\Gamma$ is an ellipse with foci $w=\pm1$ and the mapping $f$ is
$$f(z)=\sin\left(\frac{\pi{\bf F}(\frac{z}{r},r^2)}{2{\bf K}(r^2)}\right),$$ where ${\bf F}(z,k)$ is the
first elliptic integral, $${\bf F}(z,k) = \int_0^z \frac{dq}{\sqrt{(1- q^2)(1-k^2 q^2)}},$$ and ${\bf
K}(k)={\bf F}(1,k)$. The eccentricity of the ellipse~$\Gamma$ equals $\lambda=1/f(1)$. The exterior
mapping is just the Joukowski mapping $$\varphi(\zeta)=\frac1{2} \left(c_\lambda \zeta + \frac{1}{c_\lambda
\zeta}\right), \qquad c_\lambda := \frac{1+ \sqrt{1 - \lambda^2}}{\lambda},$$ and $$v_0(z) =
\frac{1}{(\varphi^{-1} \circ f)\der(z)} = \frac{2r {\bf K}(r^2) \sqrt{(r^2-z^2)(z^2 - r^{-2})}}{\pi z}= \frac{2 {\bf K}(r^2)|r^2-z^2|}{\pi }.$$
\end{remark}

\section{Conformal welding for a class of circle diffeomorphisms}\label{Poly}

Consider a diffeomorphism $\gamma:\UTs\to\UTs$ such that the function $v_0:=1/(\gamma^{-1})\der$
 has the form $v_0(z)=\sum_{k=-n}^na_kz^k$, in which case,
since $v_0$ is positive, $a_{-k}=\conj{a_k}$, and so we have two equivalent representations:
\begin{equation}\label{v_0_poly}
v_0(z)=a_0+\sum_{k=1}^na_kz^k+\frac{\conj{a_k}}{z^{k}}= \varkappa\prod_{k=1}^n\frac{e^{-it_k}}z(r_ke^{it_k}-z)(z-e^{it_k}/r_k),
\end{equation}
where  $r_k\in(0,1)$, $t_k\in\Real$, $k=1,\ldots,n$, and the coefficients $\varkappa$ and $a_k$'s are subject
to  the conditions $v_0>0$ and $\int_0^{2\pi}dt/v_0(e^{it})=2\pi$.

The set of all diffeomorphisms~$\gamma$ satisfying the above condition is dense in many important  spaces
of circle homeomorphisms.  Let us consider the problem of finding the function~$f\in\clSt$ corresponding
to $v_0$ given by~\eqref{v_0_poly}. In general, for a diffeomorphism $\gamma\in\Cc$, $\alpha\in(0,1)$, the
conformal welding is given by a unique solution to the equation
\begin{equation}\label{en_gang_til}
I_f[v_0](z):=-\frac{1}{2\pi i}\int_\UT\left(\frac{sf'(s)}{f(s)}\right)^2\frac{v_0(s)}{f(s)-f(z)}\frac{ds}{s}=0,\quad z\in\UD,
\end{equation}
regarded as an equation with respect to $f\in\clSH$. The existence and uniqueness of the solution
to~\eqref{en_gang_til} is implied by Theorem~\ref{thrm_main} and the fact that for any $\gamma\in\QS$
there exists a unique conformal welding with $f\in\clS$.

If $v_0$ is of the form~\eqref{v_0_poly}, then~\eqref{en_gang_til} can be substantially simplified by
means of calculus of residues. The residue of the expression under the integral at $s=z$ equals
$zf'(z)v_0(z)/(f(z))^2$ and the residue at the origin is of the form $P_0\big(1/f(z)\big)/f(z)$, where
$P_0$ is a polynomial of degree~$n$ with coefficients depending on~$a_k$'s and the first Taylor
coefficients of~$f$. It follows that the function $w=f(z)$ satisfies the differential equation
\begin{equation}\label{diff_eq}
 \frac{w^{n-1}dw}{P(w)}=\frac{z^{n-1}dz}{Q(z)},
\end{equation}
where $$P(w):=b_0\prod\limits_{k=1}^{n}(w-w_k),$$ $b_0$ and $w_k$'s are unknown parameters and
$$Q(z):=z^nv_0(z)=\varkappa\prod_{k=1}^n\frac{|z_k|}{z_k}(z_k-z)(z-1/\conj{z_k}),\quad z_k:=r_ke^{it_k}.$$
Since $f$ is univalent and analytic in $\UD$, $w_k$'s are exactly the images of $z_k$'s and we can suppose
that they are numbered so that $w_k=f(z_k)$.

For simplicity we suppose that all the roots of~$Q$ are simple. Then $w_k\neq w_j$ for $k\neq j$ and
comparing residues of $z^{n-1}/Q(z)$ and $f'(z)\big(f(z)\big)^{n-1}/P(f(z))$ we obtain the following
system of algebraic equations:
\begin{equation}\label{w_k}
 \frac{w_k^{n-1}}{P_k(w_k)}=A_k,\quad k=1,\ldots,n,
\end{equation}
where $$P_k(w):=\frac{P(w)}{w-w_k},\quad A_k:=\res_{z=z_k}\frac{z^{n-1}}{Q(z)}.$$ Using the residue
theorem we further conclude that $$ \frac1{b_0}=\sum_{k=1}^{n}A_k=\int_0^{2\pi}\frac{dt}{v_0(e^{it})}=1.
$$ In view of~\eqref{diff_eq} the condition $f'(0)=1$ results in the equality
\begin{equation}\label{w_k_norm}
\prod_{k=1}^nw_k=(-1)^n Q(0)=\varkappa\prod_{k=1}^n\frac{z_k}{|z_k|}.
\end{equation}

Now we can summarize the above consideration as following
\begin{proposition}
Suppose $\gamma\in\DiffS$ is such that $v_0:=1/(\gamma^{-1})\der$ is of the form~\eqref{v_0_poly}. Then
the function~$f\in\clSt$ that corresponds to $\gamma$ via conformal welding, is a solution to differential
equation~\eqref{diff_eq} with ${b_0:=1}$ and ${w_k:=f(z_k)}$. Moreover,  the vector $(w_1,\ldots,w_n)$
satisfies system~\eqref{w_k},\,\eqref{w_k_norm}, provided all the roots $z_k$ of $Q$ are simple.
\end{proposition}

\begin{remark}
Given any non-vanishing values of the parameters $w_k$, ${k=1,\ldots,n}$, differential
equation~\eqref{diff_eq} with $b_0:=1$ has a unique analytic solution ${w=w(z)}$ in a neighborhood
of~$z=0$ that satisfies the condition~${w(0)=w'(0)-1=0}$. At the same time the number of solutions of
system~\eqref{w_k},\,\eqref{w_k_norm} grows drastically as $n$ increases.
\end{remark}

The simplest case $n=1$ corresponds to the subgroup $\Mob\subset \DiffS$ consisting of  M\"obius transformations of the
unit disk restricted to $\UTs$ (excluding rotations, which correspond to $n=0$) and $f$ has the form
$z/(1-c_1 z)$, $|c_1|\in(0,1)$.  But even for $N=2$ the expressions turn out to be quite complicated.

\noindent
{\it Address:}  \newline {\sc Department of Mathematics \newline 
University of Bergen \newline Johannes Brunsgate 12, \newline Bergen 5008, \newline Norway}

\vspace{1cm}
\noindent
{\it E-mails:}
\newline
{\sc Erlend Grong:}  Erlend.Grong@math.uib.no
\newline
{\sc Pavel Gumenyuk:}  Pavel.Gumenyuk@math.uib.no
\newline
{\sc Alexander Vasil'ev:}   Alexander.Vasiliev@math.uib.no

\end{document}